\documentclass[10pt]{article}
\usepackage[mathscr]{eucal}

\usepackage{graphicx} 
\usepackage{url}
\usepackage{stmaryrd}

\newcommand\DMO[2]{\DeclareMathOperator{#1}{#2}}

\usepackage{amsmath,amsthm,amsfonts,amscd,amssymb}
\usepackage[usenames,dvipsnames]{color}

\usepackage[all]{xy} 

\usepackage{makeidx}

\newcounter{enumitemp}

\numberwithin{equation}{section}

\newcommand\ds\displaystyle

\theoremstyle{plain}

\newtheorem*{HyperbolicityTheorems}{Hyperbolicity Theorems}

\newtheorem*{theorem*}{Theorem}
\newtheorem*{loxtheorem}{Loxodromic Classification Theorems for $\FS(\Gamma;\A)$}
\newtheorem*{loxconjecture}{Loxodromic Classification Conjectures for $\CFFS(\Gamma;\A)$}

\newtheorem*{proposition*}{Proposition}

\newtheorem*{corollary*}{Corollary}

\newtheorem*{conjecture*}{Conjecture}

\theoremstyle{definition}

\newtheorem{problem}{Problem}
\newtheorem{question}[problem]{Question}

\newcommand\TOAT{\emph{Two Over All Theorem}}
\newcommand\STOAT{\emph{Strong Two Over All Theorem}}

\DeclareMathOperator{\Out}{Out}
\DeclareMathOperator{\Aut}{Aut}
\DeclareMathOperator\Inn{Inn}

\DeclareMathOperator{\Stab}{Stab}
\DeclareMathOperator\diam{diam}

\DeclareMathOperator{\MCG}{\mathsf{MCG}}
\DeclareMathOperator\Teich{Teich}
\DeclareMathOperator\Isom{Isom}

\newcommand\hyp{{\mathbf H}}
\newcommand\Z{{\mathbf Z}}

\newcommand\G{{\mathcal G}}

\newcommand{\bdy}{\partial}
\newcommand{\from}{\colon}

\newcommand\union{\cup}

\newcommand\abs[1]{\left| #1 \right|}

\newcommand\intersect{\cap}

\newcommand\subgroup{<}

\newcommand\C{\mathcal C}

\renewcommand\L{\mathcal L}

\newcommand\T{\mathcal T}

\newcommand\A{\mathscr A}
\newcommand\B{\mathscr B}
\newcommand\F{{\mathscr F}}

\newcommand\Fell[1]{\F_{\!ell} #1}
\newcommand\FellS{{\Fell  S}}
\newcommand\FellT{{\Fell T}}
\newcommand\FellU{{\Fell U}}

\newcommand\CFFS{\mathcal{F\!F}}

\newcommand\FFC{{\mathcal{F}}}

\newcommand\<\langle
\renewcommand\>\rangle
\newcommand\wh\widehat
\newcommand\disjunion\sqcup

\newcommand\act\curvearrowright

\newcommand\X{\mathscr{X}}
\newcommand\CV\X

\newcommand\BH{\cite{BestvinaHandel:tt}}

\newcommand\BookOneTag{BFH:TitsOne}
\newcommand\BookOne{\cite{\BookOneTag}}

\newcommand\RelFSOneTag{HandelMosher:RelComplexHyp}
\newcommand\RelFSOne{\cite{\RelFSOneTag}}

\newcommand\RelFSTwoTag{HandelMosher:RelComplexHypII}
\newcommand\RelFSTwo{\cite{\RelFSTwoTag}}

\newcommand\RelFSThreeTag{HandelMosher:RelComplexHypIII}
\newcommand\RelFSThree{\cite{\RelFSThreeTag}}

\newcommand\LymanCTTag{Lyman:CT}

\newcommand\ti {\tilde}

\DMO\Core{Core}
\DMO\ACore{{\hat{\mathcal{C}}}}
\DMO\truss{truss}
\newcommand\wt\widetilde
\renewcommand\O{{\mathcal O}}

\newcommand\FS{\mathcal{F\!S}}
\newcommand\FF{\mathcal{F\!F}}
\newcommand\PFS{\mathcal{P\!F\!S}}
\newcommand\CPFS{\mathcal{C\!P\!F\!S}}

\newcommand\collapsesto\succ
\newcommand\collapse\collapsesto
\newcommand\collapses\collapsesto
\newcommand\expandsto\prec
\newcommand\expand\expandsto
\newcommand\expands\expandsto

\newcommand\Kspine{\mathscr K}
\newcommand\Ks\Kspine

\newcommand\relA{\emph{rel}~$\A$}

\title{Relative free splitting and free factor complexes: \\ An overview}

\author{Michael Handel and Lee Mosher}

\begin{document}

\maketitle

\begin{abstract}
For any group $\Gamma$ and any free factor system~$\mathscr A$ of $\Gamma$, the relative outer automorphism group $\Out(\Gamma;\mathscr A)$ acts naturally on the relative free splitting complex $\FS(\Gamma;\mathscr A)$ and on the complex of relative free factor systems $\CFFS(\Gamma;\mathscr A)$, generalizing the well known actions of $\Out(F_n)$ on the absolute free splitting complex $\FS(F_n)$ and the absolute free factor complex $\FFC(F_n)$ of the rank~$n$ free group~$F_n$.

     This overview summarizes a three part work regarding the large scale geometry of $\FS(\Gamma;\mathscr A)$ and $\CFFS(\Gamma;\mathscr A)$ and the geometric dynamics of the actions on these complexes by elements of $\Out(\Gamma;\mathscr A)$. In Part I \cite{\RelFSOneTag} we prove hyperbolicity of $\FS(\Gamma;\mathscr A)$ and of $\CFFS(\Gamma;\mathscr A)$. In Parts II and III \cite{\RelFSTwoTag,\RelFSThreeTag} we study the relation between the geometric dynamics of an element of $\Out(\Gamma;\mathscr A)$ and the dynamics of its relative train track representatives. The main tool in Part II is the \TOAT, expressing an exponential flaring property of Stallings fold paths in $\FS(\Gamma;\mathscr A)$. The main tools in Part III are \emph{filling paths}, used to formulate and prove a strong version of the \TOAT.
\end{abstract}

%
%
For any group $\Gamma$ and any free factor system~$\mathscr A$ of $\Gamma$, the relative outer automorphism group $\text{Out}(\Gamma;\mathscr A)$ acts naturally on the relative free splitting complex $\mathcal{F\!S}(\Gamma;\mathscr A)$ and on the complex of relative free factor systems $\mathcal{F\!F}(\Gamma;\mathscr A)$, generalizing the well known actions of $\text{Out}(F_n)$ on the absolute free splitting complex $\mathcal{F\!S}(F_n)$ and the absolute free factor complex ${\mathcal{F}}(F_n)$ of the rank~$n$ free group~$F_n$.

     This overview summarizes a three part work regarding the large scale geometry of $\mathcal{F\!S}(\Gamma;\mathscr A)$ and $\mathcal{F\!F}(\Gamma;\mathscr A)$ and the geometric dynamics of the actions on these complexes by elements of $\text{Out}(\Gamma;\mathscr A)$. In Part I arXiv:1407.3508 we prove hyperbolicity of $\mathcal{F\!S}(\Gamma;\mathscr A)$ and of $\mathcal{F\!F}(\Gamma;\mathscr A)$. In Parts II and III \cite{\RelFSTwoTag,\RelFSThreeTag} we study the relation between the geometric dynamics of an element of $\text{Out}(\Gamma;\mathscr A)$ and the dynamics of its relative train track representatives. The main tool in Part II is the \emph{Two Over All Theorem}, expressing an exponential flaring property of Stallings fold paths in $\mathcal{F\!S}(\Gamma;\mathscr A)$. The main tools in Part III are \emph{filling paths}, used to formulate and prove a strong version of the \emph{Two Over All Theorem}.

Consider a free factorization of a group $\Gamma$ of the form $\Gamma = A_1 * \cdots * A_k * B$ ($k \ge 0$), where each $A_i \subgroup \Gamma$ is nontrivial $(1 \le i \le k)$ and the subgroup $B \subgroup \Gamma$, called the \emph{cofactor}, is free of finite rank, possibly trivial. The set of conjugacy classes $\A = \{[A_i] \,:\, 1 \le i \le k\}$ is called a \emph{free factor system} of~$\Gamma$. The outer automorphism group $\Out(\Gamma) = \Aut(\Gamma) / \Inn(\Gamma)$ acts naturally on the set of free factor systems. Fixing a particular free factor system~$\A$ of~$\Gamma$, the subgroup of $\Out(\Gamma)$ that preserves~$\A$ is denoted $\Out(\Gamma;\A)$ and is called the \emph{outer automorphism group of $\Gamma$} \relA. The foundational special case $\Out(F_n)$ arises when $k=0$ hence $\A=\emptyset$, and when $\Gamma=B \approx F_n$ is free of some rank~$n \ge 2$. 

Culler and Vogtmann, in an analogy with the Teichm\"uller space $\Teich(S)$ of a surface $S$ on which its mapping class group $\MCG(S)$ acts properly discontinuously, constructed the \emph{outer space} on which $\Out(F_n)$ acts properly discontinuously, together with an invariant spine on which the restricted action of $\Out(F_n)$ is cocompact \cite{CullerVogtmann:moduli}. This construction was subsequently generalized to a relative outer space $\O(\Gamma;\A)$ on which the group $\Out(\Gamma;\A)$ acts: first in special cases by Collins \cite{Collins:symmetric} and by McCullough and Miller \cite{McCulloughMiller:symmetric} (see below for further details of these cases); and then in general by Guirardel and Levitt \cite{GuirardelLevitt:outer}.
For a still broader generalization to deformation spaces of group actions on simplicial trees (which we shall not consider), see \cite{Forester:Deformation,GuirardelLevitt:DefSpaces}.

The curve complex $\C(S)$ of a surface $S$, on which $\MCG(S)$ acts naturally, was introduced by Harvey to serve as a kind of ``space at infinity'' associated to $\Teich(S)$ \cite{Harvey:Boundary}, in analogy to Tits buildings associated to symmetric spaces. Masur and Minsky \cite{MasurMinsky:complex1} proved hyperbolicity of $\C(S)$ and they analyzed the geometric dynamics of elements of the mapping class group $\MCG(S)$ acting on $\C(S)$: pseudo-Anosov mapping classes act loxodromically; every other mapping class has a fixed point. These results led to further studies of the large scale geometry of mapping class groups and of the geometric dynamics of their elements and subgroups, some of which explored quantitative aspects of geometric dynamics: 
\cite{GHKL:Lipschitz,BestvinaFujiwara:bounded,Bowditch:tight,Mangahas:UniformUniform,Mangahas:recipe,GHKL:Lipschitz,Fujiwara:Subgroups}. In particular, Bowditch found a positive lower bound to the stable translation lengths of the actions of pseudo-Anosov elements of $\MCG(S)$ on $\C(S)$ \cite{Bowditch:tight}. Also, see below for a further discussion (and question) regarding the work of Fujiwara \cite{Fujiwara:Subgroups}.

Certain analogues of surface curve complexes, on which the group $\Out(F_n)$ acts naturally, were proved to be hyperbolic, including: the free factor complex $\FFC(F_n)$~\cite{BestvinaFeighn:FFCHyp}; the free splitting complex $\FS(F_n)$~\cite{HandelMosher:FreeSplittingHyperbolic}; and the cyclic splitting complex \cite{Mann:CyclicSplittingComplex}. Coarse Lipschitz relations amongst such complexes were explored in \cite{KapovichRafi:HypImpliesHyp}. Analyses of geometric dynamics for actions on $\FFC(F_n)$ by elements of $\Out(F_n)$ were carried out in \cite{BestvinaFeighn:FFCHyp,Mangahas:UniformUniform}, and for actions on $\FS(F_n)$ in \cite{HandelMosher:FreeSplittingLox}. These results have provided a basis for further studies of $\Out(F_n)$ regarding large scale geometry and geometric dynamics of elements and of subgroups, often by analogy with $\MCG(S)$. Here is a small selection of such studies, with a highlight on ones that involve quantitative conclusions: \cite{BestvinaFeighn:subfactor,Taylor:subfactor,ClayUlyanik:simultaneous,Webb:HypNotCatZero,KMPT:RandomOut,HandelMosher:BddCohomology,BCH:Connectivity}. We expect that the results for $\Out(\Gamma;\A)$ presented here --- regarding hyperbolicity of $\FS(\Gamma;\A)$ and $\CFFS(\Gamma;\A$), and quantititative geometric dynamics of the actions on these spaces by elements of $\Out(\Gamma;\A)$  --- will provide a similar basis for further studies.

\tableofcontents

\section{Summary of the main results} Parts I, II, III of this series are organized around versions of the theorems for $\Out(F_n)$ mentioned above, generalized to the setting of relative outer automorphism groups $\Out(\Gamma;\A)$. 

Part~I proves hyperbolicity of $\FS(\Gamma;\A)$ known as the ``free splitting complex of $\Gamma$ rel~$\A$'', and of~$\CFFS(\Gamma;\A)$ known as the ``complex of free factor systems of $\Gamma$ rel~$\A$'' or more briefly as the ``complex of relative free factor systems''.

Parts~II and~III describe the geometric dynamics of actions on $\FS(\Gamma;\A)$ by individual elements of $\Out(\Gamma;\A)$, with some partial results for actions on $\CFFS(\Gamma;\A)$ included in the final section of Part~III. Parts~II and~III are distinguished from each other by their methods of proof, as alluded to in their subtitles: the \TOAT\ in Part II; and concepts of \emph{filling paths} in Part III. These methods yield applications  to the action of $\Out(\Gamma;\A)$ on $\FS(\Gamma;\A)$ that are new even for the case of $\Out(F_n)$ acting on $\FS(F_n)$: see Theorems~A and~B (discussed at length below) for more refined conclusions regarding quantitative geometric dynamics; and see Part~II Section 3 (with a brief discussion in Section 4 below) for the coarse Lipschitz property of the natural ``inclusion'' map into $\FS(\Gamma;\A)$ defined on the relative outer space $\O(\Gamma;\A)$.

\smallskip

We turn to a brief summary of basic concepts needed to understand $\FS(\Gamma;\A)$ and $\CFFS(\Gamma;\A)$ and their actions by $\Out(\Gamma;\A)$, with further details found in Part I. After that we break into separate descriptions of the theorems of Parts I, II and~III.

The inclusion partial order on subgroups of $\Gamma$ induces a partial order on the set of free factor systems of~$\Gamma$, denoted $\A \sqsubset \B$. One can apply Grushko's Theorem and the Kurosh Subgroup Theorem together to deduce that if $\Gamma$ is finitely generated then this partial ordering has a unique minimum~$\A$ (see Part I, Section 2.4), in which case $\Out(\Gamma;\A) = \Out(\Gamma)$; this is the special case covered in \cite{McCulloughMiller:symmetric}. A different special case, covered in \cite{Collins:symmetric}, is for $F_n = \<a_1,\ldots,a_n\>$ and $\A = \bigl\{[\<a_1\>],\ldots,[\<a_n\>]\bigr\}$, where $\Out(F_n;\A)$ is known as the group of \emph{symmetric} outer automorphisms of $F_n$.

Fixing a free factor system~$\A$ of $\Gamma$, the partial order can be restricted to free factor systems $\B$ \relA\ --- meaning $\A \sqsubset \B$ --- and then further restricted by requiring that $\A \ne \B \ne \{[\Gamma]\}$. The simplicial realization of this poset is the \emph{complex of relative free factor systems} denoted $\CFFS(\Gamma;\A)$. The natural action of $\Aut(\Gamma)$ on the set of free factors of $\Gamma$ induces actions of $\Out(\Gamma;\A)$ first on the set of free factor systems \relA, and then on $\CFFS(\Gamma;\A)$. Embedded in $\CFFS(\Gamma;\A)$ is the ``(relative) free factor complex'' $\mathcal F(\Gamma;\A) \subset \CFFS(\Gamma;\A)$, which is obtained by further restricting the partial order to those free factor systems rel~$\A$ of the form $\B = \A' \union \{[B]\}$ where $\A' \sqsubset \A$, and $[B] \not\in \A$; in this situation we refer to $B$ as a \emph{free factor rel~$\A$}. In all but the simplest cases the inclusion map $\mathcal F(\Gamma;\A) \hookrightarrow \CFFS(\Gamma;\A)$ is an $\Out(\Gamma;\A)$ equivariant quasi-isometry (see Part~I \cite[Proposition 6.3]{\RelFSOneTag}). We can thus think of $\CFFS(\Gamma;\A)$ as an indirect generalization of the free factor complex $\mathcal F(F_n)$ on which $\Out(F_n)$ acts. For our purposes it is generally easier to work with free factor \emph{systems} rel~$\A$ than with individual free factors rel~$\A$, but we will occasionally need to move back and forth between those two points of view.

A \emph{free splitting} of $\Gamma$ \relA\ is a minimal action $\Gamma \act T$ on a simplicial tree $T$ with finitely many edge orbits, with trivial edge stabilizers, and with a constraint on the stabilizer subgroup $\Stab(v)$ of each vertex $v$ of $T$:  either $\Stab(v)$ is trivial; or there exists a subgroup $A \subgroup \Stab(v)$ such that $[A] \in \A$. For any free splitting $T$ rel~$\A$, the set of conjugacy classes of nontrivial vertex stabilizers in $T$ forms a free factor system of $\Gamma$ \relA\ denoted $\FellT$, called the \emph{elliptic free factor system} of~$T$. If $\FellT=\A$ then one says that $T$ is a Grushko free splitting \relA. Two free splittings are equivalent when there is a $\Gamma$-equivariant piecewise linear homeomorphism between them. The \emph{collapse} partial order $[S] \collapsesto [T]$ on set of equivalence classes of free splittings \relA\ is defined by existence of a \emph{collapse map} $S \mapsto T$ which is a continuous, \hbox{$\Gamma$-equivariant,} piecewise linear map $S \mapsto T$ such that the inverse image of each point is connected. The simplicial realization of this partial ordering is the \emph{relative free splitting complex} $\FS(\Gamma;\A)$. For any free splitting $\Gamma \act T$ \relA\ and any $\phi \in \Out(\Gamma;\A)$, the action of $\Gamma$ on~$T$ may be precomposed by any automorphism $\Phi \in \Aut(\Gamma;\A)$ that representes $\phi \in \Out(\Gamma;\A)$, obtaining another free splitting \relA\ that depends only on $\phi$ up to equivalence; this induces a natural simplicial action of $\Out(\Gamma;\A)$ on~$\FS(\Gamma;\A)$.  There is a natural, $\Out(\Gamma;\A)$ equivariant ``inclusion'' $\O(\Gamma;\A) \hookrightarrow \FS(\Gamma;\A)$ of the relative outer space that identifies $\O(\Gamma;\A)$ with the dense open set complementary to the nowhere dense subcomplex of $\FS(\Gamma;\A)$ consisting of all non-Grushko free splittings~\relA.

\smallskip

The hyperbolicity theorems of Part~I \cite{\RelFSOneTag} give two different analogues of the hyperbolicity of curve complexes: one for $\FS(\Gamma;\A)$, generalizing hyperbolicity of $\FS(F_n)$ \cite{HandelMosher:FreeSplittingHyperbolic}; and one for $\CFFS(\Gamma;\A)$, generalizing hyperbolicity of $\mathcal F(F_n)$ \cite{BestvinaFeighn:FFCHyp}. Also, the proofs in   Part~I use Stallings fold paths to describe a family of quasigeodesics, generalizing what was done for $\FS(F_n)$ in \cite{HandelMosher:FreeSplittingHyperbolic}, and for $\FFC(F_n)$ in \cite{KapovichRafi:HypImpliesHyp}. 

\begin{HyperbolicityTheorems} For any group $\Gamma$ and free factor system~$\A$ of $\Gamma$ the following hold:
\begin{enumerate}
\item\label{ItemFSTheorem} The relative free splitting complex $\FS(\Gamma;\A)$ is hyperbolic. The set of Stallings fold paths in $\FS(\Gamma;\A)$ forms a coarsely transitive family of reparameterized quasigeodesics in $\FS(\Gamma;\A)$. 
\item\label{ItemFFTheorem}
The complex of relative free factor systems $\CFFS(\Gamma;\A)$ is hyperbolic (outside of certain low complexity exceptions). Under a natural $\Out(\Gamma;\A)$-equivariant, Lipschitz projection map $\FS(\Gamma;\A) \to \CFFS(\Gamma;\A)$, the set of Stallings fold paths in $\FS(\Gamma;\A)$, projects to a coarsely transitive family of reparameterized quasigeodesics in $\CFFS(\Gamma;\A)$.
\end{enumerate}
\end{HyperbolicityTheorems}
\noindent
A ``reparameterized'' quasigeodesic is a quasigeodesic precomposed by a monotonic surjection from one subinterval of the number line to another; see \cite[Section~5.1]{\RelFSOneTag}. ``Coarse transitivity'' of a family of paths in a space means that for any two points $x,y$ in that space, there is a path in the family with initial point near $x$ and terminal point near~$y$.

In the \emph{Hyperbolicity Theorems} above there are several implicit constants, and there will be several explicit constants in Theorems~A and~$B$ below. Most constants in our results will depend in a \emph{particular manner} only on $\Gamma$ and~$\A$, expressed as $K = K(\Gamma,\A)$, with the meaning that $K$ depends only on two specific numerical invariants of $\Gamma$ and~$\A$, namely: the cardinality $k \ge 0$ of the free factor system \hbox{$\A = \{[A_i] \,:\, 1 \le i \le k\}$}; and the rank of any cofactor of~$\A$ (although cofactors of $\A$ are not unique, not even up to conjugacy in $\Gamma$, nonetheless they are well defined ``modulo~$\A$'' \cite[Lemma~2.2]{\RelFSOneTag}, and so the rank of a cofactor is well-defined). Constants of the form $K(\Gamma;\A)$ are therefore independent of the isomorphism types of the free factors $A_i$ $(1 \le i \le k)$. For~example, the implicit constants in the \emph{Hyperbolicity Theorems} are all of this form: hyperbolicity constants; quasigeodesic constants; coarse transitivity constants; and so on. Also, the constants in Theorem~A and~B are of this form: stable translation length bounds; and so on. 

\smallskip

In Parts~II and~III we analyze the ``geometric dynamics'' of elements of $\Out(\Gamma;\A)$ acting on the hyperbolic space $\FS(\Gamma;\A)$ (and see below for a limited discussion of the action on $\FF(\Gamma;\A)$). Every self-isometry $\phi \from X \to X$ of a metric space has a \emph{stable translation length} $\tau_\phi = \lim_{n \to \infty} \frac{1}{n} d(x,\phi^n(x)) \ge 0$, well-defined independent of $x \in X$. If $X$ is Gromov hyperbolic then $\tau_\phi > 0$ if and only if $\phi$ is \emph{loxodromic} meaning that for each $x \in X$ the orbit map $n \mapsto \phi^n(x)$ is a quasi-isometric embedding $\Z \mapsto X$. The nonloxodromic case $\tau_\phi = 0$ may be refined into two subcases: $\phi$ is \emph{elliptic} if every $\phi$-orbit is bounded; otherwise $\phi$ is \emph{parabolic}. 
 
The main results of Parts~II and~III relate relative train track dynamics to quantitative geometric dynamics for the actions of each $\phi \in \Out(\Gamma;\A)$ on $\FS(\Gamma;\A)$ and on $\CFFS(\Gamma;\A)$, in particular to bounds on stable translation lengths and on orbit diameters of~$\phi$. Relative train track representatives for $\Out(F_n)$ were introduced by Bestvina and Handel in \BH, with further developments by Bestvina, Feighn and Handel in \BookOne\ including attracting laminations and their expansion factors, and with later extensions to $\Out(\Gamma;\A)$ by Lyman \cite{Lyman:RTT,Lyman:CT}. These theories arose in analogy to the dynamical theory of surface mapping class groups $\MCG(S)$, introduced by W.\ Thurston \cite{FLPKM,CassonBleiler} where pseudo-Anosov mapping classes are characterized as having an invariant measured geodesic lamination that fills the surface and that attracts all nearby measured geodesic laminations. 

\begin{loxtheorem}[Part~II \cite{HandelMosher:RelComplexHypII} and Part~III \cite{HandelMosher:RelComplexHypIII}] \quad\hfill

\vskip -10pt
\noindent
There exist constants $A=A(\Gamma;\A) > 0$, $B=B(\Gamma;\A) > 0$ and $\Omega=\Omega(\Gamma;\A) \ge 1$ such that for any $\phi \in \Out(\Gamma;\A)$ the following hold:

\smallskip\noindent\textbf{Theorem A:} If $\phi$ has a filling attracting lamination $\Lambda$ with expansion factor $\lambda_\phi>1$ then 
$$A \le \tau_\phi \le B \, \log(\lambda_\phi)
$$
\textbf{Theorem B:} The following are equivalent:
\begin{enumerate}
\item $\phi$ has a filling attracting lamination.
\item $\phi$ acts loxodromically on $\FS(\Gamma;\A)$ (equivalently $\tau_\phi>0$).
\item $\phi$ acts with unbounded orbits on $\FS(\Gamma;\A)$.
\item\label{ItemOrbitsUniformlyLarge}
Every orbit of the action of $\phi$ on $\FS(\Gamma;\A)$ has diameter $\ge \Omega$.
\end{enumerate}
\end{loxtheorem} 
\noindent
These two theorems taken together may be regarded as giving a quantitative refinement of the following qualitative statement, which generalizes the main result of \cite{HandelMosher:FreeSplittingLox} for $\Out(F_n)$: 
\begin{corollary*}
$\phi \in \Out(\Gamma;\A)$ acts loxodromically on $\FS(\Gamma;\A)$ if and only if $\phi$ has a filling lamination; otherwise $\phi$ acts elliptically. 
\end{corollary*}
\noindent
This corollary evidently follows from Theorem~B, although we \emph{prove} the implication (1)$\implies$(2) of Theorem~B --- the ``if'' direction of the corollary --- by direct application of the lower bound of Theorem~A (see the \emph{Outline of the proof of Theorem~B} in the introduction of Part III). 

\medskip

In the final section of Part III we formulate analogues of Theorems A and~B for the actions on $\CFFS(\Gamma;\A)$ by elements of $\Out(\Gamma;\A)$; these may be regarded as quantitative refinements of a theorem of Guirardel and Horbez saying that $\phi \in \Out(\Gamma;\A)$ acts loxodromically on $\CFFS(\Gamma;\A)$ if and only if $\phi$ is fully irreducible rel~$\A$, \emph{assuming} that certain sporadic cases are ruled out (for the full statement see Theorem 5.1 of Part III \cite{HandelMosher:RelComplexHypIII}). But we can prove only certain partial results regarding these analogues; the rest is left as questions and conjectures.

\smallskip

In what follows of this \emph{Overview} we first give some history of Parts I, II and~III, next we discuss the mathematical methods used for the proofs in each part, and finally we discuss our original motivation for the \TOAT, arising from still incomplete investigations of the large scale geometry of the \emph{automorphism} group $\Aut(F_n)$ and of its relativizations $\Aut(\Gamma;\A)$.

\section{Some history of this work}

\textbf{A brief history of Part I.} As mentioned earlier, the hyperbolicity theorems for $\FS(F_n)$ found in \cite{HandelMosher:FreeSplittingHyperbolic}, and for $\FFC(F_n)$ found in \cite{BestvinaFeighn:FFCHyp}, were motivated by hyperbolicity of the curve complex $\C(S)$ associated to each finite type surface~$S$, found  in~\cite{MasurMinsky:complex1}. The latter theorem played a key element in the development by Masur and Minsky of a ``hierarchy theory'' for the mapping class group $\MCG(S)$ \cite{MasurMinsky:complex2}. That theory as a whole cracked open the large scale geometry of the mapping class group $\MCG(S)$, by applying hyperbolicity of the curve complexes of $S$ itself and of connected essential subsurfaces of $S$, together with various relations amongst those complexes. As an example of such a relation, given a connected essential subsurface $\Sigma \subset S$, the curve complex of $\Sigma$ serves as a kind of ``natural complement'' of the curve complexes of the components of $S-\Sigma$; for example, the product of those curve complexes embeds naturally as a subcomplex of the whole curve complex.

Pursuing further analogies with hierarchy theory has proved fruitful in studying the large scale geometry of $\Out(F_n)$, for example applications of hyperbolicity of $\FS(A)$ and $\FFC(A)$ for \emph{free factors} $A \subgroup F_n$ as found in \cite{BestvinaFeighn:subfactor}. But will this be enough for further developments? 

We came to suspect that for any free factor system $\A = \{[A_1],\ldots,[A_K]\}$ of~$F_n$, hyperbolic complexes associated to the individual free factors $A_1,\ldots,A_K$, such as $\FS(A_1),\ldots,\FS(A_K)$, should have a ``natural complement'' for large scale geometry purposes, and that purpose might be served by some kind of ``relative'' complex. Motivated by the work of Guirardel and Levitt cited earlier on relative outer spaces, together with the observation that the relative free splitting complex $\FS(F_n;\A)$ is the simplicial completion of the relative outer space $\O(F_n;\A)$ --- we were led to consider $\FS(F_n;\A)$ for the role of a natural complement (whether it \emph{does} fill that role is beyond what we can say right now). We thus set out to generalize our proof of hyperbolicity of $\FS(F_n)$ found in \cite{HandelMosher:FreeSplittingHyperbolic} so as to work for $\FS(F_n;\A)$, soon realizing that the proof should also work more generally for $\FS(\Gamma;\A)$. The result was Part I of this work \RelFSOne, first posted on the arXiv in 2014 (see also \cite{Horbez:HyperbolicGraphs}).

\smallskip\textbf{A brief history of Parts II and III.} We originally intended to follow up Part I with a generalization of the work \cite{HandelMosher:FreeSplittingLox}, obtaining a classification of the \emph{qualitative} dynamics of $\phi \in \Out(\Gamma;\A)$ acting on $\FS(\Gamma;\A)$: $\phi$ acts loxodromically if and only if $\phi$ has a filling lamination \relA, otherwise $\phi$ acts elliptically. Meanwhile, while working on the methods that eventually led to the \TOAT\ and the \STOAT\ (with other motivations discussed in later sections of this \emph{Overview}), we came to realize that these methods might have applications to the \emph{quantitative} dynamics of $\phi \in \Out(\Gamma;\A)$ acting on $\FS(\Gamma;\A)$, obtaining explicit bounds to various constants that are only implicit in qualitative dynamics. We slowly formulated the pieces of Theorems~A and~B, at first just for the case of $\Out(F_n)$. But as we came to realize that such applications should work in general for $\Out(\Gamma;\A)$, we realized that it would make no sense to write about quantitative dynamics only for the special case of $\Out(F_n)$. This led to a long delay after first posting Part~I \RelFSOne, while we slowly developed the mathematics that resulted in Parts II and III with applications to Theorems~A and~B.

\section{Methods of Part I: Hyperbolicity \cite{\RelFSOneTag}} 

\paragraph{Hyperbolicity of $\FS(\Gamma;\A)$.} This proof, in \cite[Sections 4 and 5]{\RelFSOneTag}, follows the proof of hyperbolicity of $\FS(F_n)$ in \cite{HandelMosher:FreeSplittingHyperbolic}, but with changes taken up from \cite{BestvinaFeighn:subfactor} as discussed in the introduction of \RelFSOne\ under the heading \emph{Comparison with earlier methods of proof}; see also \cite{Horbez:HyperbolicGraphs}. 

Here, outside of very brief comments on the application of the Masur--Minsky hyperbolicity axioms in Part~I, we focus on a description of those tools that are used for the proofs in Part~I and are also used in Parts~II and~III of this work, namely: \emph{fold paths}; and diagrammatic tools involving fold paths called \emph{projection diagrams} and \emph{component free splitting units}. 

\smallskip
\emph{Fold paths.} Fold paths are our primary tools for analyzing the geometry of $\FS(\Gamma;\A)$. They were introduced by Stallings \cite{Stallings:folding,Stallings:folds} to solve problems in group theory. A general theory of fold paths is described in~\cite{BestvinaFeighn:bounding}; for those elements of that theory relevant to our current setting, see \cite[Section 4.1]{\RelFSOneTag}, including a form of Stallings' fold theorem \cite[Theorem~4.3]{\RelFSOneTag}. 

A ~\emph{foldable map} $f \from S \to T$ is a $\Gamma$-equivariant map between two free splittings (of $\Gamma$ \relA) such that $S$ is covered by open arcs on which $f$ restricts to an injection. For any such $f$, consider a factorization of $f$ as a sequence of $\Gamma$-equivariant maps between free splittings as follows: 
$$\xymatrix{
S = T_I \ar[r]_{f_{I+1}} \ar@/^2pc/[rrr]^{f} & T_{I+1} \ar[r]_{f_{I+2}} & \cdots \ar[r]_{f_{K}} &T_K = T
}$$
The ``foldable'' property is inherited by every factor $f^i_j = f_j \circ \cdots \circ f_{i+1} \from T_i \to T_j$ ($I \le i < j \le K)$ along the whole sequence 
\cite[Lemma~4.2~(3)]{\RelFSOneTag}.
Stallings fold theorem ensures the existence of such a factorization which is so fine that each $f_i \from T_{i-1} \to T_i$ is a \emph{fold map}, determined by folding together some pair of short oriented closed arcs in $T_{i-1}$ having the same initial point; the result is a \emph{fold sequence}, also called a \emph{fold factorization} of the given foldable map $f \from S \to T$. The source and target of any fold map represent a pair of vertices of $\FS(\Gamma;\A)$ that are connected by an edge path of length at most~$2$, and so one may visualize a fold sequence as a path in \hbox{the $1$-skeleton} of $\FS(\Gamma;\A)$, hence the terminology of a \emph{fold~path}. 

\smallskip
\emph{The Masur--Minsky axioms.} These axioms for verifying hyperbolicity of a given connected simplicial complex, found in \cite{MasurMinsky:complex1}, require two further mathematical objects to be given. First, one must give a \emph{coarsely transitive} family of paths, meaning a family which (up to bounded error) connects any two vertices; for this purpose we use Stallings fold paths in $\FS(\Gamma;\A)$. One must also give a ``projection function'' for each path in the given family, mapping the whole complex to that path; we define these functions using ``projection diagrams'' (\cite[Section 5.2]{\RelFSTwoTag}, reviewed briefly just below). The axioms themselves describe geometric properties of the family of paths and their projection functions. The main theorem of \cite{MasurMinsky:complex1} lets one conclude that the given complex is hyperbolic, once the axioms are verified. For $\FS(\Gamma;\A)$ we carry out these verifications in \cite[Section 5]{\RelFSTwoTag}. From that same theorem it also follows that fold paths are reparameterized quasigeodesics; but we also describe a particularly concrete and useful quasigeodesic reparameterization using ``free splitting units''.

\smallskip
\emph{Combing rectangles.} Defining projection maps and free splitting units, and verifying the Masur--Minsky axioms, are accomplished using certain commutative diagrams of maps between free splittings, which one may visualize as subcomplexes of the 2-skeleton of $\FS(\Gamma;\A)$. The fundamental building block is a \emph{combing rectangle} (see \cite[Section 4.3]{\RelFSOneTag}), and when two such rectangles are stacked in a particular way one obtains a \emph{collapse--expand rectangle}:
$$
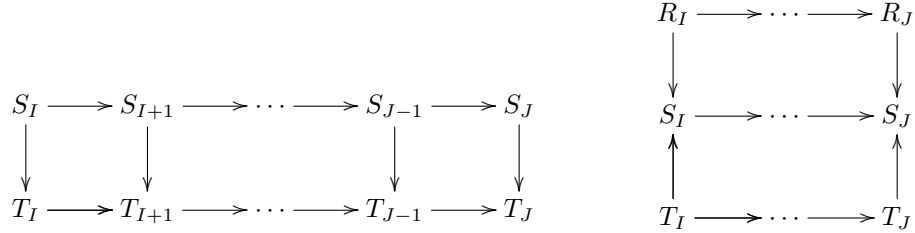
\xymatrix{
 \\
S_I 
	\ar[r]
	\ar[d]
 & S_{I+1} 
 	\ar[d]
	\ar[r]
 & \cdots 
 	\ar[r]
 & S_{J-1} 
 	\ar[d]
	\ar[r]
 & S_J 
 	\ar[d]
	\\
T_I 
	\ar[r]
 	\ar[r]
& T_{I+1}                                   
	\ar[r]
& \cdots 
	\ar[r]
 & T_{J-1}                                 
 	\ar[r]
 & T_J
}
\qquad\qquad
\xymatrix{
R_I \ar[r] \ar[d]
                 & \cdots \ar[r]           
                 & R_J  \ar[d]
                                                                           \\
S_I \ar[r]                         
                 & \cdots \ar[r]           
                 & S_J & \\
T_I 
	\ar[r]
	\ar[u]
                 	\ar[r]
                 	\ar[u]
                 & \cdots 
                          \ar[r]
                 & T_J \ar[u]
}
$$
The left diagram depicts a combing rectangle: each horizontal row is a foldable sequence; each vertical arrow is a collapse map; and a compatibility relation is required amongst the collapsed edges of the free splittings~$S_i$ along the top row. In a context where the bottom row and right side are given, the combing rectangle is constructed in \cite[Lemma~4.9]{\RelFSOneTag} and is called a \emph{combing-by-expansion} rectangle. And when instead the top row and the right side are given (and the diagram is then usually flipped across a horizontal axis so that vertical arrows point upward), the combing rectangle is constructed in \cite[Lemma~4.8]{\RelFSOneTag} and is called a \emph{combing-by-collapse} rectangle. A \emph{collapse--expand} rectangle (on the right) is a stack of two combing rectangles, with a combing--by--collapse rectangle on the bottom and a combing--by--expansion rectangle on~top.

\smallskip
\emph{Projection diagrams.} Given any fold path $T_I \mapsto\dots\mapsto T_K$ and any free splitting $R$ representing a vertex of $\FS(\Gamma;\A)$, the projection function from $\FS(\Gamma;\A)$ to that fold path is evaluated on $R$ with the aid of \emph{projection diagrams} from $R$ to $T_I \mapsto\dots\mapsto T_K$ \cite[Definition~5.1]{\RelFSOneTag}, each of which has the form 
$$
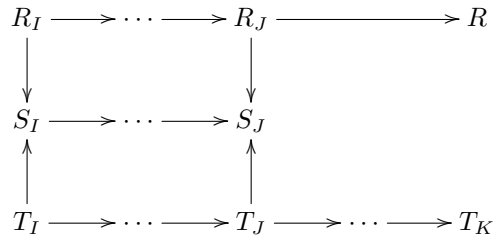
\xymatrix{
R_I \ar[r] \ar[d] & \cdots \ar[r] & R_{J} \ar[rr] \ar[d] &&  R \\
S_I \ar[r]          & \cdots \ar[r] & S_{J} \\
T_I \ar[r] \ar[u] & \cdots \ar[r] & T_{J} \ar[r] \ar[u]  & \cdots \ar[r] & T_K \\
}$$
In this diagram, the bottom row is the given fold sequence, the middle and top rows are foldable sequences, the right endpoint of the top row is the given free splitting $R$, and the portion of the diagram between columns $I$ and $J$ is a collapse--expand rectangle. By visualizing this diagram in the 2-skeleton of $\FS(\Gamma;\A)$ one can see a combinatorial expression of a fellow traveller property between the given fold path $T_I \mapsto \cdots \mapsto T_K$ and a certain fold path which starts near $T_I$ --- namely, at $R_I$ --- and which ends at $R$. The value of $J$ expresses where one loses track of the fellow travelling property between the two fold paths. Poorly chosen combinatorics can ignore geometry: some projection diagrams may have (much) larger values of $J$ than others. The~projection of $R$ to $T_I \mapsto \cdots\mapsto T_K$ is defined by an optimization: it is that position $T_\Delta$ such that $\Delta \in \{I,\ldots,K\}$ maximizes $J$ over all projection diagrams from $R$ to $T_I \mapsto\cdots\mapsto T_K$. 

At~the heart of the proof of hyperbolicity, found in Part I \cite[Section 5.5]{\RelFSOneTag}, is a diagrammatic argument that starts with an optimal projection diagram from $R$ to $T_I \mapsto\cdots\mapsto T_K$ together with a finite geodesic segment in the $1$-skeleton of $\FS(\Gamma;\A)$ having one endpoint at $R$ and having opposite endpoint at some other vertex~$R'$. The argument proceeds by making efficient use of this data to produce a projection diagram from $R'$ to $T_I \mapsto\cdots\mapsto T_K$. This is done by constructing, analyzing, and cutting and pasting, a ``big diagram'' of combing rectangles that is built atop the foldable sequence $R_I \mapsto\cdots \mapsto R_J \mapsto R$ by stacking collapse--expand rectangles whose right sides are arrows along the given geodesic from $R$ to $R'$ (see \cite[Section 5.5, Figure 3]{\RelFSOneTag}).

\smallskip
\emph{Free splitting units.} Along a Stallings fold path in $\FS(F_n)$, free splitting units were developed in \cite{HandelMosher:FreeSplittingHyperbolic} as an explicit, combinatorially defined quasigeodesic reparameterization. We do the same here in $\FS(\Gamma;\A)$, together with an improved concept of ``component free splitting units'' that is simpler in expression and more easily applied. For basic definitions and for the proof that component free splitting units give a quasigeodesic reparameterization of Stallings fold paths, see Part I \cite[Definition 4.17, Lemma 4.18, and Corollary 5.5]{\RelFSOneTag}. Free splitting units are applied in two places in this work: the proof of hyperbolicity of $\FS(\Gamma;\A)$ in Part I; and the proof of Theorem~B in Part~III \cite[Section 3.6]{\RelFSThreeTag}. 

For a brief account of component free splitting units \cite[Definition 4.7]{\RelFSOneTag}, consider a fold path \hbox{$T_I \mapsto\cdots\mapsto T_J$.} Given~a collapse expand diagram atop that fold path with top row $R_I \mapsto\cdots\mapsto R_J$ (as depicted earlier), as $i$ varies from $I$ to $J$ we are concerned with the evolution of \emph{nondegenerate subforests} of~$R_i$, meaning $\Gamma$-invariant subforests in which no component is a single point. When a nondegenerate subforest $\beta_J \subset R_J$ is given at the end of the top row, one can pull it back along the arrows in the top row, producing a sequence of nondegenerate subforests \hbox{$\beta_i \subset R_i$ ($I \le i \le J$)} \cite[Section 4.3]{\RelFSOneTag}. To say that \emph{there~is $<1$ component free splitting unit between $T_I$ and~$T_J$} means that there is a collapse expand diagram atop the given fold path, and there exists nondegenerate subforest $\beta_J \subset R_J$ having exactly one orbit of components, such that when $\beta_J$ is pulled back along the arrows in the top row each subforest $\beta_i$ also has exactly one orbit of components. More generally, the number of component free splitting units between $T_I$ and $T_J$ is the maximum integer $\Upsilon \ge 0$ for which there exists a sequence $I \le i(0) < i(1) < \cdots < i(\Upsilon) \le J$ such that for all $1 \le u \le \Upsilon$ there is $<1$ component free splitting unit between $T_{i(u-1)}$ and $T_{i(u)}$.

\paragraph{Hyperbolicity of $\CFFS(\Gamma;\A)$.} This proof, found in \cite[Section 6]{\RelFSOneTag}, uses a method developed by Kapovich and Rafi for ``pushing forward'' hyperbolicity under a particular class of maps. They applied this method to a natural Lipschitz map $\pi \from \FS(F_n) \mapsto \mathcal F(F_n)$, using it to deduce hyperbolicity of $\mathcal F(F_n)$ as a consequence of hyperbolicity of $\FS(F_n)$~\cite{KapovichRafi:HypImpliesHyp}. Applying the same method, we use a natural Lipschitz map $\pi \from \FS(\Gamma;\A) \to \CFFS(\Gamma;\A)$ to deduce hyperbolicity of $\CFFS(\Gamma;\A)$ from hyperbolicity of $\FS(\Gamma;\A)$, in the non-exceptional cases where $\CFFS(\Gamma;\A)$ is of dimension~$\ge 1$ (see \cite[Proposition 6.2]{\RelFSOneTag}; in certain very low complexity exceptional cases the complex $\CFFS(\Gamma;\A)$ can be $0$-dimensional and infinite, hence disconnected). 

The Lipschitz map $\pi$ inputs a vertex of $\FS(\Gamma;\A)$ represented by a free splitting $T$ of $\Gamma$ \relA. In the special case where $\FellT \ne \A$ one simply takes $\pi[T]=\FellT$; this guarantees that $\pi$ maps the vertex set of $\FS(\Gamma;\A)$ surjectively to the vertex set of $\CFFS(\Gamma;\A)$. In the opposite case where $\FellT = \A$, equivalently $T$ is a Grusko free splitting rel~$\A$, one considers the finitely many free splittings $U$ such that $\FellU \ne \A$ and such that there exists a collapse map $T \mapsto U$. The set of all such free factor systems $\FellU$ is nonempty, and is collected together to form a subset $\Pi[T] \subset \CFFS(\Gamma;\A)$. In the case where $\FellT = \A$, one then chooses $\pi[T]$ arbitrarily from the subset $\Pi[T]$. The fact that $\pi$ is coarsely well-defined and Lipschitz follows from \cite[Lemma 6.4]{\RelFSOneTag} which gives a uniform bound $\diam\bigl(\Pi[T]\bigr) \le 4$; the desired path of length~$\le 4$ needed to verify this bound is explicitly constructed in a case-by-case manner.

The key hypothesis of the Kapovich--Rafi method that must be verified says that for any pair of free splittings $S,T$ representing vertices $[S],[T] \in \FS(\Gamma;\A)$, if their projections $\pi[S],\pi[T] \in \CFFS(\Gamma;\A)$ have distance~$\le 1$ then for any $\FS(\Gamma;\A)$ geodesic path connecting $[S]$ and $[T]$ the $\pi$-image of that path has uniformly bounded diameter in $\CFFS(\Gamma;\A)$. One then exploits the collapse maps that arise in the ``opposite case'' of the definition of $\pi$, in order to reduce the proof to the special case that $\A \ne \FellS \sqsubset \FellT$. One may then replace a given geodesic path from $S$ to $T$ with a Stallings fold path $S = T_0 \mapsto \cdots \mapsto T_K=T$, a replacement that is justified by the fact that the latter path is a reparameterized quasigeodesic. Combining $\A \ne \FellS = \FellT_0 \sqsubset \FellT_K$ with $\Gamma$-equivariance it follows that $\A \ne \FellT_0 \sqsubset \cdots \sqsubset \FellT_K$, and so the vertices $\{\pi[T_0],\ldots,\pi[T_K]\} = \{\FellT_0,\ldots,\FellT_K\}$ span a single simplex of $\CFFS(\Gamma;\A)$ which has diameter~$\le 1$.

\section{Methods of Part II: Stable translation lengths and \hfill\break the Two Over All Theorem \cite{\RelFSTwoTag}} 

The central result of Part II \RelFSTwo\ is the \TOAT, which describes an exponential flaring property for a foldable map between two free splittings: 

\begin{theorem*}[The \TOAT.] There is a constant $\Delta=\Delta(\Gamma;\A)$ such that for any foldable map $f \from S \to T$ between two free splittings of $\Gamma$ \relA, if the distance between $S$ and $T$ in $\FS(\Gamma;\A)$ is at least $n \Delta$ then the free splitting $S$ has two natural edges $E_1,E_2$, in two different orbits of the action of $\Gamma$, such that for each natural edge $E \subset T$ and each $i=1,2$ the image path $f(E_i)$ crosses $2^{n-1}$ distinct natural edges of $T$ in the orbit of $E$. 
\end{theorem*}

For the story of our original motivation for this theorem, see Section~\ref{SectionMotivationForTOAT}, the questions therein, and in particular the final heading \emph{``Flaring properties and the \TOAT''}. Meanwhile, although the questions in Section~\ref{SectionMotivationForTOAT} remain open, the \TOAT\ and its strengthened version discussed below have some interesting applications, which gives us the opportunity of presenting these theorems in this work.

\paragraph{Applications.} We give two applications of the \TOAT\ in Part II \cite{\RelFSTwoTag}. One application, in Section~4, is a proof of the upper bound on stable translation lengths in the conclusion of Theorem~A. The other application, in Section~3, is concerned with the natural topological embedding $\O(\Gamma;\A) \hookrightarrow \FS(\Gamma;\A)$ of the relative outer space into the free splitting complex. The image of this embedding is an open dense subset consisting of all Grushko free splittings, and the complement of the image is a subcomplex of $\FS(\Gamma;\A)$. The conclusion in this setting is that the embedding is Lipschitz, with respect to the log-Lipschitz semi-metric on $\O(\Gamma;\A)$ \cite{FrancavigliaMartino:TrainTracks,Meinert:TrainTrackMaps} and the usual simplicial metric on $\FS(\Gamma;\A)$. As a quick corollary, the natural systole map $\O(\Gamma;\A) \mapsto \CFFS(\Gamma;\A)$ is also coarse Lipschitz, generalizing the classical case $\O(F_n) \mapsto \FFC(F_n)$ which follows from work of Bestvina and Feighn \cite{BestvinaFeighn:FFCHyp}.

In both applications, the idea of the proof is simply to express the desired conclusion as an exponential growth property. The upper bound $B \log(\lambda_\phi)$ in the application to Theorem~A expresses an exponential growth property of a relative train track representative of $\phi \in \Out(\Gamma;\A)$ that witnesses the value of the expansion factor $\lambda_\phi$. The Lipschitz property in the second application expresses an exponential growth property of an optimal map that witnesses the value of the semi-metric. In both cases, the proof comes down to working through the definitions in order to match the exponential growth conclusion of the \TOAT\ with the corresponding exponential growth expression of the desired application.

\paragraph{Proof of the \TOAT.} The proof starts with a simple reduction argument, employing a Stallings fold factorization of $f$ to reduce the theorem to its ``uniterated'' form, meaning the special case where $n=1$. In that case, the final conclusion is replaced by the statement that each of the paths $f(E_i)$ $(i=1,2)$ crosses some natural edge in every natural edge orbit of~$T$. That uniterated version was first worked out by us for a special case formulated in the setting of $\FS(F_n)$, namely for a foldable map between two Grushko free splittings --- in other words, a lift to universal covers of a foldable map between two marked graphs. We provide a sketch of this special version in the first pages of Section~5 of \RelFSTwo. While that sketch gives many of the key ideas needed for proving the \TOAT\ in the general context of $\FS(\Gamma;\A)$, it nonetheless took us quite a bit more time even to just formulate the theorem in its general context, let alone to prove it.

For now we mention just one key strategy of the proof, namely to exploit the fact that Stallings fold factorizations of a given foldable map are usually far from unique. This failure of uniqueness is an impediment in certain applications; see for example \cite{Skora:deformations,GuirardelLevitt:outer} which contain proofs of contractibility of outer spaces that work by making a very careful choice of Stallings fold paths depending continuously on their endpoints. But this failure of uniqueness is a boon for our application: given a foldable map $f \from S \to T$, by sifting through the many choices one finds opportunities for doing lots of folds without moving far in $\FS(\Gamma;\A)$. For example, suppose that $\beta \subset S$ is a proper $\Gamma$-invariant subforest. In this situation one can define a ``partial'' Stallings fold factorization $f \from S = T_0 \mapsto\cdots\mapsto T_k \mapsto T$ by \emph{prioritizing} folds of $\beta$, meaning that a fold factor $T_{i-1} \mapsto T_i$ is allowed to fold only edges of $T_{i-1}$ that lie in image of $\beta$. When the priority process stops at $T_K$ and no more such folds are possible for the image $\beta' \subset T_K$ of $\beta$, it follows that the set of vertices of $\FS(\Gamma;\A)$ represented by the free splittings $\{T_0,\ldots,T_k\}$ has diameter~$\le 2$ (see \cite[Proposition 5.4]{\RelFSTwoTag}). Furthermore, it follows that the map $T_k \mapsto T$ is injective on each component of $\beta'$, and that injectivity is itself a useful property that one can exploit as the proof proceeds. The proof of the \TOAT\ essentially consists of a carefully chosen sequence of fold priority steps, each satisfying a uniform diameter bound, such that the number of steps in the sequence also satisfies a uniform bound.

\section{Methods of Part III: Stable translation lengths and \hfill\break filling paths \cite{\RelFSThreeTag}}

The central theme of Part III is the study of filling \emph{laminations} using the concept of filling \emph{paths}, with applications to the proof of the lower bound of Theorem~A and to the proof of Theorem~B. 

The set of attracting laminations of an outer automorphism $\phi \in \Out(\Gamma;\A)$ can be concretely described using any choice of relative train track representative of $\phi$. In the foundational case of $\Out(F_n)$, such representatives were described in \cite{BestvinaHandel:tt} as maps on marked graphs. In the case of $\Out(\Gamma;\A)$, the description was generalized in \cite{Lyman:CT} as maps on marked graphs-of-groups. One can also lift to Bass-Serre trees (i.e.\ to universal covers in the case of $\Out(F_n)$), describing a relative train track representative as a map \hbox{$F \from T \to T$} defined on a Grushko free splittings $T$ of $\Gamma$ \relA, where $F$ is required to be ``$\Phi$-twisted equivariant'' with respect to a choice of automorphism $\Phi \in \Aut(\Gamma;\A)$ representing~$\phi$ (see \cite[Sections 2.1 and 4.3]{\RelFSTwoTag}. For simplicity in this overview we make the special assumption that $F$ is ``EG-aperiodic'', meaning that its transition matrix has a positive power, and implying that $\phi$ has a unique attracting lamination~$\Lambda$ (in practice one reduces to this special case using constructions found in \cite[Section 4.3.3]{\RelFSTwoTag}). The generic leaves of $\Lambda$ are those lines $L$ in the tree $T$ that are \emph{exhausted by iteration tiles of~$F$}, meaning that $L$ is the union of some nested sequence of paths of the form $F^k(E_k)$ ($k \ge 1$) where each $E_k \subset T$ is an~edge.  For any other free splitting $R$ (of $\Gamma$ \relA) one can use the $\Gamma$-equivariant relation between the end spaces of $T$ and $R$ to realize each leaf of $\Lambda$ in $R$ \cite[Lemma 3.1]{\RelFSThreeTag}. Given an  attracting lamination $\Lambda$ of $\phi$, to say that $\Lambda$ \emph{fills} $\Gamma$ \relA\ is equivalent to saying that in every free splitting $R$ of $\Gamma$ \relA, the union of the set of lines in $R$ that realize leaves of $\Lambda$ is the entire tree~$R$ \cite[Corollary 3.5]{\RelFSThreeTag}.  

For the implication~(4)$\implies$(1) of Theorem~B, after first using (4) to identify an appropriate attracting lamination $\Lambda$ of $\phi$, one needs to \emph{prove} that $\Lambda$ fills $\Gamma$ \relA\ \cite[Section 3]{\RelFSThreeTag}. And for the lower bound of Theorem~A (from which follows the implication (1)$\implies$(2) of Theorem~B), after first finding an appropriate value $A = A(\Gamma;\A)>0$, one needs to \emph{apply} that $\Lambda$ fills $\Gamma$ \relA\ as part of the proof of the inequality $\tau_\phi \ge A$ \cite[Section 4]{\RelFSThreeTag}. 

For purposes of both \emph{proving} and \emph{applying} the property that $\Lambda$ fills $\Gamma$ \relA, our main tools are \emph{filling paths} in free splittings. Given a free splitting $S$ of $\Gamma$ \relA\ and a finite, embedded path $\alpha \subset S$, to say that $\alpha$ \emph{fills} $S$ means that for any collapse map $f \from R \to S$ defined on some free splitting $R$ of $\Gamma$ \relA, letting $\tilde\alpha \subset R$ denote the unique minimal path in $R$ such that $f(\tilde\alpha)=\alpha$, the following holds: for every natural edge $E \subset R$ there exists $\gamma \in \Gamma$ such that the interior of the path $\tilde\alpha$ contains the natural edge $\gamma \cdot E$ (see \cite[Definition 2.2]{\RelFSThreeTag}). 

The theory of filling paths is developed in \cite[Section 2]{\RelFSThreeTag}, leading up to the following theorem, which builds on the \TOAT\ by strengthening its conclusions at the expense of using a larger constant:

\begin{theorem*}[The \STOAT] There is a constant $\Theta=\Theta(\Gamma;\A)$ such that for any foldable map $f \from S \to T$ between two free splittings of $\Gamma$ \relA\ and for any $n \ge 1$, if the distance between $S$ and $T$ in $\FS(\Gamma;\A)$ is at least $n\Theta$ then the free splitting $S$ has two natural edges $E_1,E_2$, in two different orbits of the action of $\Gamma$, such that each image path $f(E_i)$ has $2^{n-1}$ non-overlapping subpaths each of which fills~$T$.
\end{theorem*}
\noindent
The proof of this theorem is found in \cite[Section 2.4]{\RelFSThreeTag}, with an outline given in the introduction to \RelFSThree. The proof combines an application of the original \TOAT\ with results from the theory of filling paths \cite[Sections 2.1--2.3]{\RelFSThreeTag}.

\paragraph{Proving the implication (4)$\implies$(1) of Theorem~B.} This proof is found in Section 3 of Part III, with an outline given in the introduction to Part III. Here is a very brief account which emphasizes the applications of the \STOAT, of free splitting units, and of projection diagrams.

The goal is to find an appropriate value of $\Omega=\Omega(\Gamma;\A) > 0$ and use it to prove that for all $\phi \in \Out(\Gamma;\A)$, if every orbit of $\phi$ in $\FS(\Gamma;\A)$ has diameter $\ge \Omega$ then $\phi$ has a filling attracting lamination. By a straightforward relative train track argument, a preliminary lower bound $\Omega_1 = 5$ on orbit diameters implies the existence of an attracting lamination $\Lambda$ of $\phi$, an EG-aperiodic train track representative $F \from T \to T$ of $\phi$ which is associated to $\Lambda$ as described above ($T$ is possibly non-Grushko), and a ``suspension axis'' of $F$ which is a certain bi-infinite, $\phi$-periodic Stallings fold path in which the given train track representative $F$ appears as a ``first return map'' (see \cite[Section 3.5]{\RelFSThreeTag}). We may write this axis and its fold maps in the form
$$\cdots \xrightarrow{f_{-K}} \underbrace{T_{-K}}_{T \cdot \phi^{-1}} \xrightarrow{f_{-K+1}} T_{-K+1} \xrightarrow{f_{-K+2}} \cdots T_{-1} \xrightarrow{f_0} \underbrace{T_0}_T \xrightarrow{f_1} T_1 \xrightarrow{f_2} \cdots \xrightarrow{f_{K-1}} T_{K-1} \xrightarrow{f_K} \underbrace{T_K}_{T \cdot \phi} \xrightarrow{f_{K+1}} \cdots
$$
In this notation, the ``first return map'' on $T=T_0$ is the composition of the foldable map $f_K \circ \cdots \circ f_1 \from T = T_0 \to \cdots \to T_K = T \cdot \phi$ with the twisted equivariant PL homeomorphism $T \cdot \phi \to T$. By using the constant of the \STOAT\ as a second lower bound $\Omega_2 = \Theta$ to orbit diameters, one proves that this axis of $\phi$ is a bi-infinite quasigeodesic in the hyperbolic space $\FS(\Gamma;\A)$ --- note that this proves conclusion (2) of Theorem~B, as a step along the way towards proving conclusion (1). Component free splitting units play an important role in this deduction: with an appropriate choice of the lower bound on orbit diameters applied to a pair of points on that axis, we obtain a lower bound on component free splitting units, which leads to uniform increase of distance along that axis. Finally, given an arbitrary free splitting $R$ of $\Gamma$ \relA, one \emph{very carefully} chooses a finite subsegment of the suspension axis of~$F$, encompassing a sufficiently broad region surrounding a closest point projection from $R$ to that quasigeodesic (see the \emph{Restricted Projection Property} and the subheading \emph{Further Remarks} in \cite[Section 3.7]{\RelFSThreeTag}). A~projection diagram from $R$ to that subsegment can then be used, in conjunction with a bounded cancellation argument, to prove that $R$ is covered by the set of lines in $R$ that realize leaves of $\Lambda$. This proves that $\Lambda$ fills $\Gamma$ \relA, assuming the lower bound $\Omega=\max\{\Omega_1,\Omega_2\}$ on diameters of orbits of~$\phi$.

\paragraph{Proving the lower bound of Theorem~A.} This proof is found in Section 4 of Part III, also with an outline in the introduction to Part III. The goal is to find an appropriate constant $A=A(\Gamma;\A)>0$ and use it to prove that if $\phi \in \Out(\Gamma;\A)$ has a filling lamination $\Lambda$ then the translation length $\tau_\phi$ of its action on $\FS(\Gamma;\A)$ is $\ge A$.

The proof starts with arguments from train track theory that produce a uniform exponent $\kappa_0=\kappa_0(\Gamma;\A)$ such that for each $\phi$ as above there exists an EG-aperiodic (possibly non-Grushko) train track representative $F \from T \to T$ having the following properties: generic leaves of $\Lambda$ are exhausted by iteration tiles of $F$; and all entries of the transition matrix of $F^{\kappa_0}$ are~$\ge 4$, in other words for any two edges $E,E' \subset T_0$ the tile $F^{\kappa_0}(E)$ crosses at least $4$ different edges in the orbit of $E'$. This exponent $\kappa_0=\kappa_0(\Gamma;\A)$ is referred to as a ``PF-exponent'' of $\Gamma$ \relA\ (see \cite[Definition 3.7 and Proposition 4.3]{\RelFSThreeTag}). 

The hard work of the proof is to apply the hypothesis that $\Lambda$ fills $\Gamma$ \relA, together with elements of the theory of filling paths, in order to find a uniform integer valued factor $\mu=\mu(\Gamma;\A)$ such that $\omega_0 = \mu \kappa_0$ is a \emph{filling exponent}; this means (roughly speaking) that with notations as above, for each edge $E \subset T_0$, the tile $F^{\omega_0}(E)$ fills $T_0$ (see \cite[Definition 3.10 and Proposition 4.4]{\RelFSThreeTag}). The underlying idea of this proof is to consider one such edge $E$, and to keep track of the sequence of ``filling ranks'' of iterated tiles $F^{m \kappa_0}(E)$. These ranks must form an increasing sequence, and that sequence must increase strictly until the moment that it becomes constant. But once it does become constant, one uses the assumption that $\Lambda$ fills to prove that $F^{m\kappa_0}(E)$ is a filling path.

\paragraph{Geometric dynamics on $\CFFS(\Gamma;\A)$: Conjectures and partial results.} In the final Section~5 of Part~III, we turn to questions of geometric dynamics on the complex of relative free factor systems $\CFFS(\Gamma;\A)$, acted on by elements of $\Out(\Gamma;\A)$. The \emph{qualitative} dynamics are given in the following theorem of Guirardel and Horbez, generalizing the foundational case for $\FFC(F_n)$ proved by Bestvina and Feighn in \cite[Theorem 9.3]{BestvinaFeighn:FFCHyp}, and the extension to $\FFC(F_n;\A)$ by Gupta in \cite{Gupta:Loxodromic}. 

To state the result, one defines $\phi \in \Out(\Gamma;\A)$ to be \emph{fully irreducible} \relA\ when the only $\phi$-periodic free factor systems \relA\ are $\A$ itself and $\{[\Gamma]\}$. For various reasons we rule out a small number of ``exceptional'' cases (as defined in Section 2.5 of Part I).

\begin{theorem*}[\protect{\cite[Section 4]{GuirardelHorbez:SubgroupClassification}}]
Assuming that $(\Gamma;\A)$ is nonexceptional, for each $\phi \in \Out(\Gamma;\A)$, the action of $\phi$ on $\mathcal F(\Gamma;\A)$ is loxodromic if and only if $\phi$ is fully irreducible \relA, and otherwise its action is elliptic. \qed
\end{theorem*}
\noindent
For discussion of the few ``exceptional'' cases see Section 5 of~Part~III.

Regarding \emph{quantative} dynamics on $\mathcal F(\Gamma;\A)$, we conjecture the following refinement of the above theorem, which is formulated as an analogue, in the context of $\CFFS(\Gamma;\A)$, of Theorems~A and~B regarding $\FS(\Gamma;\A)$. 

\begin{loxconjecture} There exist constants $A^\FF > 0$, $B^\FF > 0$, and $\Omega^\FF > 0$, each depending only on $\Gamma$ and $\A$, such that the following hold:
\begin{description}
\item [Conjecture A${}^\FF$:] If $\phi$ is fully irreducible with expansion factor $\lambda_\phi > 1$ then
$$A^\FF \le \tau^\FF_\phi \le B^\FF \log(\lambda_\phi)
$$
\item[Conjecture B${}^\FF$:] The following are equivalent (the equivalence (1)$\iff$(2) being known already as cited above):
\begin{enumerate}
\item $\phi$ is fully irreducible 
\item $\phi$ acts loxodromically on $\CFFS(\Gamma;\A)$; equivalently, $\tau^\FF_\phi > 0$.
\item $\phi$ acts with unbounded orbits on $\CFFS(\Gamma;\A)$.
\item Every orbit of the action of $\phi$ on $\CFFS(\Gamma;\A)$ has diameter $\ge \Omega^\FF$.
\end{enumerate}
\end{description}
\end{loxconjecture}
\noindent
We briefly describe some partial results towards this conjecture. 

As happened for the case of $\FS(\Gamma;\A)$, the lower bound of Conjecture A${}^\FF$ may be regarded as strengthening the (known) implication (1)$\implies$(2); also, the implications (2)$\implies$(3)$\implies$(4) of Conjecture B${}^\FF$ are immediate from the definitions. 

\smallskip
In \cite[Section 5]{\RelFSThreeTag} we shall consider the remaining pieces of the above conjectures: the~upper bound in Conjecture~A${}^\FF$ has a reasonably quick proof by combining the corresponding upper bound from Theorem~A with a Lipschitz projection argument; and we shall prove the implication (4)$\implies$(1) for the special case $\Gamma=F_n$ using a relative train track argument. And although some initial steps of the latter argument extend to general groups $\Gamma$ using the relative train track theory of Lyman \cite{\LymanCTTag}, that theory as it stands is not (yet) strong enough to extend the full argument; for a brief discussion of this issue, in Part III see the footnote to the paragraph following Reduction~5.7. 

To summarize, the unresolved portions of Conjectures~A${}^\FF$ and B${}^\FF$ are:
\begin{itemize}
\item The lower bound in Conjecture~A${}^\FF$, 
\item The implication (4)$\implies$(1) for general groups $\Gamma$. 
\end{itemize}
We propose the following question as a possible avenue of attack on the lower bound (see \cite[Section 5]{\RelFSThreeTag} for further discussion):

\begin{question}
Does there exist a combinatorial measurement of ``free factor units'' along a Stallings fold path, in analogy to the use of free splitting units in the proof of the lower bound of Theorem~A, which gives a quasi-distance formula for the projection of that fold path to~$\FF(\Gamma;\A)$?
\end{question}

\smallskip\textbf{Further questions about stable translation lengths.} As alluded to earlier, Bowditch's study of the action of $\MCG(S)$ on the curve complex $\C(S)$ yields lower bounds for stable translation lengths of pseudo-Anosov elements, and those bounds take the quite striking form of a positive rational number whose denominator depends only on the topology of~$S$. In some \emph{Remarks} found at the end of the introduction to Part III we discuss these issues in further detail, leading to the following questions, each of which would result in similar rational lower bounds to stable translation lengths:

\begin{question}
Does there exist an exponent $m=m(\Gamma;\A)$ such that if $\phi \in \Out(\Gamma;\A)$ acts loxodromically on $\FS(\Gamma;\A)$ then $\phi^m$ has a bi-infinite \emph{geodesic} axis in $\FS(\Gamma;\A)$? 
\end{question}

\begin{question}
Does there exist an exponent $n=n(\Gamma;\A)$ such that if $\phi \in \Out(\Gamma;\A)$ acts loxodromically on $\FF(\Gamma;\A)$ then $\phi^n$ has a bi-infinite \emph{geodesic} axis in $\CFFS(\Gamma;\A)$?
\end{question}

Fujiwara studied the following question in \cite{Fujiwara:Subgroups}. Consider a group action $G \act X$ on a hyperbolic space $X$. Given two loxodromic elements $\phi,\psi \in G$ whose axes have disjoint endpoint pairs in the Gromov boundary $\bdy X$, one can show by a standard ping-pong argument that there exists a constant $M \ge 1$ such that for all $m,n \ge M$ the group elements $\phi^m,\psi^n$ freely generate a quasi-isometrically embedded subgroup of~$G$. In this general context it is interesting to ask: Under what conditions can the constant $M$ be chosen independently of $\phi$ and $\psi$? Fujiwara proved that this can be done if the action $G \act X$ is acylindrical, and so in particular when $G \act X$ is a word hyperbolic group acting on its Cayley graph, or $\MCG(S)$ acting on the curve complex $\C(S)$. 
\begin{question}\label{QuestionFujiwara}
For the action $\Out(\Gamma;\A) \act \FS(\Gamma;\A)$, can the constant $M$ be chosen independently of $\phi$ and $\psi$?
\end{question}
\noindent
We note here that the action $\Out(\Gamma;\A) \act \FS(\Gamma;\A)$ is generally \emph{not} acylindrical, and so Fujiwara's theorem does not apply. Nonetheless many \emph{consequences} of acylindricality do hold for this action, such as a positive lower bound to translation lengths (Theorem B). More such consequences are known in the restricted setting of $\Out(F_n) \act \FS(F_n)$. To explain this, consider a general hyperbolic action $G \act X$ and an element $\gamma \in G$ acting loxodromically with fixed point pair $\{P,Q\} \in \bdy X$. If the action $G \act X$ is acylindrical then $\gamma$ satisfies the WPD property; and although this fails in the setting of the action $\Out(F_n) \act \FS(F_n)$ (see \cite[Corollary 1.2]{HandelMosher:FreeSplittingLox}), on the other hand one can often verify a weaker version of WPD known as WWPD (see \cite[Theorem B]{HandelMosher:BddCohomology}). Also, if the action $G \act X$ is acylindrical then the subgroup of $G$ that jointly stabilizes $P$ and $Q$ is virtually cyclic; this again fails in the setting of $\Out(F_n) \act \FS(F_n)$, but on the other hand that subgroup can be explicitly described, as shown in \cite[Section 5]{HandelMosher:FreeSplittingLox}. Taken together, this feels like an accumulation of evidence that Question~\ref{QuestionFujiwara} is ripe for an attack.

\section{Motivating the \TOAT: \\Hyperbolicity questions for~$\Aut(F_n)$} 
\label{SectionMotivationForTOAT}
Our original motivation for the \TOAT\ arose from the following question. Given the successes at finding natural, useful hyperbolic spaces on which the \emph{outer} automorphism group $\Out(F_n)$ acts, one might wonder:

\begin{question}
\label{QuestionAutHyp}
Does there exist a natural and useful hyperbolic space on which the \emph{automorphism} group $\Aut(F_n)$ acts? 
\end{question}

There is an interesting literature on negative results regarding this question, i.e.\ proofs of \emph{non}-hyperbolicity of certain natural candidates for a positive result: see \cite{BBW:FreeFactorGraph}; and see also \cite{Hamenstadt:SpottedDisc} for the closely related ``spotted handlebody'' automorphism group. 

In~work following \cite{HandelMosher:FreeSplittingHyperbolic} we pondered a candidate for a positive answer to Question~\ref{QuestionAutHyp}, but this similarly led to a negative result: as we shall explain momentarily, the \emph{pointed free splitting complex} of $F_n$ is not hyperbolic. By examining the mode of failure, we then formulated another candidate for a positive answer to Question~\ref{QuestionAutHyp}, a space we shall refer to as the \emph{coned, pointed free splitting complex of~$F_n$}.

We then pursued some steps of a still incomplete outline for a proof of hyperbolicity of the coned, pointed free splitting complex complex, using ``flaring'' methods. This led us to formulate a flaring statement along Stallings fold paths in $\FS(F_n)$ which, generalized to $\FS(\Gamma;\A)$, became the \TOAT. We shall sketch out some details of this pursuit, first displaying a diagram of deformation spaces that motivates the definition of pointed free splitting complex. Next we explain the failure of hyperbolicity of the pointed free splitting complex, due to the presence of ``geometric axes''. We then state Questions~6 and~7 which ask whether the ``coned pointed free splitting complex'' --- obtained from the pointed free splitting complex by coning off all geometric axes --- is a candidate for an affirmative answer for Question~5. Finally, we discuss the relation between hyperbolicity and flaring methods, using that discussion to motivate the \TOAT.

\paragraph{A diagram of deformation spaces.} The diagram below shows four pairs of deformation spaces associated to two groups $\Gamma$: two pairs associated to $\Gamma = \pi_1 S$ where $S$ is a closed oriented surface $S$ of genus~$\ge 2$; and two associated to the free group~$\Gamma = F_n$ of rank~$n \ge 2$. The group $\Aut(\Gamma)$ acts on the first space of each pair, and $\Out(\Gamma)$ acts on the second space, and there is a ``forgetful'' function from the first space to the second that is equivariant with respect to the canonical homomorphism $\Aut(\Gamma) \mapsto \Out(\Gamma)$ with kernel $\Inn(\Gamma)$. The deformation spaces in the first line are metric spaces with proper actions, and in the second they are (putatively) hyperbolic spaces with improper actions.

\medskip

\begin{tabular}{c||c|c}
$\vphantom{\biggl|}$
             &$\Aut(\pi_1 S ) \mapsto \Out(\pi_1 S)$ 
             		& $\Aut(F_n) \mapsto \Out(F_n)$
								\\
								\hline\hline
$\vphantom{\Big|}$ 
	proper actions	& (Bers)				   & (Culler Vogtmann) \\
$\vphantom{\Big|}$
  				& $\T(S,p) \mapsto \T(S)$     & $\A(F_n) \mapsto \X(F_n)$ 
								\\ 
								\hline 
$\vphantom{\Big|}$
  	improper hyperbolic(?) actions	& (Kent Leininger Schleimer) &  \\
  	 & $\C(S,p) \mapsto \C(S)$ & $\fbox{?} \mapsto \FS(F_n)$ \end{tabular} 

\medskip
Bers defined the ``forgetful'' fibration $\T(S,p) \mapsto \T(S)$ of Teichm\"uller spaces \cite{Bers:FiberSpaces}, under which a conformal structure on $S-p$ with removable singularity at $p$ is extended naturally to a conformal structure on $S$. One can describe $\T(S)$ as the deformation space of faithful, discrete actions $\pi_1 S \act \hyp^2$ modulo conjugation by elements of $\Isom_+(\hyp^2) \approx \text{PSL}(2,\mathbb R)$. The space $\T(S,p)$ has a similar description but each action $\pi_1 S \act \hyp^2$ comes with a choice of basepoint $\ti p \in \hyp^2$, and the conjugating element must take basepoint to basepoint. The fibration $\T(S,p) \mapsto \T(S)$ is defined by simply forgetting the basepoint $\ti p$, so fibers are naturally identified with $\hyp^2$.

Kent, Leininger and Schleimer describe an analogous ``forgetful'' function $\C(S,p) \mapsto \C(S)$ of curve complexes: a curve system on $S-p$, upon forgetting $p$, includes into~$S$ to form a curve system on~$S$ \cite{KentLeiningerSchleimer:TreesMCG}. Given a $0$-cell of $\C(S)$ represented by a curve system $C \subset S$, one can also represent that $0$-cell by the graph of groups $\G$ that is dual to~$C$, and then by the Bass-Serre tree $T$ of $\G$ on which $F_n$ acts; we denote that $0$-cell by~$[T]$. We can think of (some canonical choice of) $T$ as the fiber itself, in the space $\C(S,p)$, over the $0$-cell $[T]$, the individual points of that fiber then being just the points of~$T$. Given a $1$-cell of $\C(S)$ represented by a nested pair of curve systems $C_0 \subset C_1$ in $S$, one can also represent that $1$-cell by a graph-of-groups collapse map $\G_1 \mapsto \G_0$, and then by a $\pi_1 S$-equivariant collapse map of Bass-Serre trees $T_1 \mapsto T_0$. In $\C(S,p)$, the preimage of that $1$-cell may be identified with the mapping cylinder of the collapse map $T_1 \mapsto T_0$, and the collapse map itself may be thought of as a ``connection map'' between the two adjacent fibers $T_1,T_0 \subset \C(S,p)$. Each of the complexes $\C(S,p)$ and $\C(S)$ is hyperbolic: the hyperbolicity proof of Masur and Minsky \cite{MasurMinsky:complex1} is agnostic about closed surfaces versus punctured surfaces. 

The Culler-Vogtmann outer space $\X(F_n)$, on which $\Out(F_n)$ acts, may be represented as the space of minimal, proper, isometric actions $\rho \from F_n \to \Isom(T)$ on metric simplicial trees~$T$, modulo conjugation of the $F_n$ action by an equivariant homothety \cite{CullerVogtmann:moduli}. The associated space on which $\Aut(F_n)$, denoted $\A(F_n)$, was christened the \emph{autre espace} by Paulin \cite{Vogtmann:OuterSpaceSurvey}. In~$\A(F_n)$, the fiber over the point $[T] \in \X(F_n)$ is (some canonical choice of)~$T$ itself, the points of the fiber being the points of~$T$.

We come now to the lower right corner of the diagram, using the free splitting complex $\FS(F_n)$ as a hyperbolic space on which $\Out(F_n)$ acts improperly. By analogy with the constructions above one obtains the \emph{pointed free splitting complex} $\PFS(F_n)$ on which $\Aut(F_n)$ acts. Over each $0$-cell of $\FS(F_n)$ represented by a free splitting action $F_n \act T$, the associated fiber of $\FS(F_n)$ is thought of as (some canonical copy of)~$T$ itself. Over each $1$-cell of $\FS(F_n)$ represented by a collapse map between two free splittings, the preimage in~$\PFS(F_n)$ can be represented as the mapping cylinder of that collapse map, and such collapse maps may be thought of as ``connection maps'' between adjacent fibers in~$\PFS(F_n)$. 

\smallskip

Unfortunately, $\PFS(F_n)$ is not a hyperbolic space.

\paragraph{Failure of hyperbolicity of $\PFS(F_n)$: From geometric axes to quasiflats.} We explain how to assemble certain quasiflats in $\PFS(F_n)$ out of certain lines in fibers of $\PFS(F_n)$ called \emph{geometric axes}, thereby proving that $\PFS(F_n)$ is not hyperbolic.

Given a marked graph $G$ and a circuit in $G$ represented as an immersion $c \from S^1 \to G$, to say that $c$ is a \emph{geometric circuit (in~$G$)} means that its mapping cylinder $S_c$ is a topological surface that supports a pseudo-Anosov homeomorphism. This mapping cylinder is of course just the quotient space
$$S_c = (S^1 \times [0,1]) \sqcup G \, \bigm/  \, (x,1) \sim c(x) \,\text{(for all $x \in S^1$)}
$$
The requirement that $S_c$ be a surface imposes a restriction on the circuit $c$: it must cross each edge of~$G$ exactly twice; and its ``Whitehead graph'' at each vertex $v$, which records the ways in which $c$ crosses a regular neighborhood of $v$, must be a circle. Also, recalling that $n \ge 2$, the pseudo-Anosov requirement imposes the restriction that if $n=2$ then $S_c$ is not a Klein bottle with one hole, hence $S_c$ is a torus with one hole. The image of the induced homomorphism \hbox{$\pi_1 S^1 \xrightarrow{c_*} \pi_1 G \approx F_n$} is a maximal infinite cyclic subgroup that is well-defined up to inner automorphism. Any such subgroup of $F_n$, associated to any choice of $G$ and $c$, is called a \emph{geometric subgroup} of $F_n$, and the conjugacy class of that subgroup is called a \emph{geometric conjugacy class} of~$F_n$. Using an argument of Stallings \cite{Stallings:transversality} one can show that a geometric subgroup \emph{fills} $F_n$ in the sense that it is not a subgroup of any proper free factor of~$F_n$. The natural inclusion $G \hookrightarrow S_c$ is a deformation retract, so there is an induced isomorphism $F_n \approx G \approx \pi_1 S_c$ well-defined up to inner automorphism of~$F_n$, thus giving $S_c$ the structure of a ``marked surface''. Equivalence of marked surfaces is defined by existence of a homeomorphism that preserves marking up to inner automorphism. Two geometric circuits $c_1$, $c_2$ in two marked graphs $G_1,G_2$ determine equivalent marked surfaces $S_{c_1}$, $S_{c_2}$ if and only if $c_1,c_2$ determine the same geometric conjugacy class of~$F_n$; this follows from the theorem of Zieschang, Vogt and Coldeway \cite{ZieschangVogtColdewey}, a relative version of the Dehn--Nielsen theorem, saying that if $h \from S_{c_1} \to S_{c_2}$ is a homotopy equivalence that restricts to a homeomorphism $\bdy S_{c_1} \to \bdy S_{c_2}$ then $h$ is homotopic rel boundary to a homeomorphism $S_{c_1} \mapsto S_{c_2}$.

Given a marked graph $G$, with its universal covering tree $T = \wt G$ regarded as the fiber in $\PFS(F_n)$ over the $0$-cell $[T] \in \FS(F_n)$, any geometric subgroup of $F_n$ acts loxodromically on the tree $T$ (because the subgroup fills $F_n$), and the axis of that action is a line $L \subset T \subset \PFS(F_n)$ called a \emph{geometric axis} in~$\PFS(F_n)$. 

\smallskip

We recall some facts about geometric outer automorphisms from \cite{BestvinaHandel:SurfaceTT}. To say that a fully irreducible element $\phi \in \Out(F_n)$ is \emph{geometric} means that there exists a marked surface $S$ and a pseudo-Anosov homeomorphism \hbox{$h \from S \to S$} that represents $\phi \in \Out(F_n) \approx \Out(\pi_1 S)$. This holds if and only if some nontrivial conjugacy class of $F_n$ is $\phi$-periodic, and if so then there is an EG-irreducible train track representative $f \from G \to G$ of (some power of) $\phi$ and a closed indivisible Nielsen path of $f$ forming a geometric circuit $c$ in~$G$, such that the restriction $f \mid c$, when straightened, yields $c$ back again; in this situation $S$ can be identified with the surface $S_c$. Up in the universal cover $T = \wt G$ the map $f$ lifts to a (twisted equivariant) map $F \from T \to T$ of $\phi$, and $c$ lifts to a geometric axis $L \subset T$, such that the restriction of $F$ to $L$, when straightened, yields $L$ back again. Consider a bi-infinite fold axis associated to $\phi$, as in the diagram following the earlier heading ``Proving the implication (1)$\implies$(4) of Theorem B''. Using notations from that diagram, each free splitting along that axis may be regarded as the fiber $T_k \subset \PFS(F_n)$ over the $0$-cell $[T_k] \in \FS(F_n)$, and there are geometric axes $L_k \subset T_k$, such that for each $k \in \Z$ the restriction of the fold map $T_{k-1} \to T_k$ to $L_{k-1}$, when straightened, yields~$L_k$. The desired quasiflat in $\PFS(F_n)$ is
$$\L = \bigcup_{k \in \Z} L_k \subset \bigcup_{k \in \Z} T_k \subset \PFS(F_n)
$$
Here is a formula for a quasi-isometric embedding $Q \from \mathbb Z^2 \to \PFS(F_n)$ whose image has finite Hausdorff distance from $\L$. As $k$ varies, the axes $L_k$ all have the same infinite cyclic stabilizer $\<\gamma\> \subgroup F_n$, where the geometric element $\gamma \in F_n$ represents the same conjugacy class represented by the circuit $c$. For each fold map $f_k \from T_{k-1} \to T_k$ we have $L_k \subset f_k(L_{k-1}) \subset N_r(L_k)$ for a uniform radius $r>0$. Using this one can extract a sequence $p_k \in L_k \subset T_k$ such that $f_k(p_{k-1})=p_k$ ($k \in \mathbb Z$). Define 
$$Q(k,l) = \gamma^l(p_k) \in T_k \subset \L
$$
The image of $Q$ clearly has finite Hausdorff distance from $\L$. Using the sup norm on $\mathbb Z^2$, clearly $Q$ is Lipschitz. What remains is to derive upper bounds on the two coordinate differences of the form
$$\abs{k-k'} \le C \, d(Q(k,l),Q(k',l')) \qquad \abs{l-l'} \le D \, d(Q(k,l),Q(k',l'))
$$
where $d(\cdot,\cdot)$ is distance in $\PFS(F_n)$. The first upper bound comes from the fact that $\phi$ acts loxodromically on $\FS(F_n)$ --- every element of $\Out(F_n)$ represented by a pseudo-Anosov homeomorphism on a marked surface with one boundary component is fully irreducible. The second upper bound can be derived from the fact that $\gamma$ fills $F_n$, and so there is a uniform positive lower bound to the translation lengths of the $\gamma$ actions on all fibers of $\PFS(F_n)$.

\paragraph{The coned, pointed free splitting complex.} The manner in which geometric axes are implicated in the failure of hyperbolicity of $\PFS(F_n)$ suggests the following question:

\begin{question}
\label{QuestionCPFS}
Is there a relatively hyperbolic structure on $\PFS(F_n)$ in which every geometric axis in every fiber has bounded diameter? 
\end{question}

For example, starting with $\PFS(F_n)$ one could attach an individual cone to every geometric axis in the fiber over every vertex of $\FS(F_n)$. Unfortunately the resulting space is still not hyperbolic, because there are ``new'' quasiflats (which we overlooked in the original draft of this overview) that one can construct by starting with one of Lyman's \emph{almost Nielsen paths} \cite{Lyman:CT}, rather than an ordinary Nielsen path as described just above. Here is a more extensive coning operation that includes all of those ``new'' quasiflats, resulting in the coned space described below in Question~\ref{QuestionRealCPFS}, which we present as a contender for a positive answer to Question~\ref{QuestionCPFS}.

Consider a free splitting $T$ and a geometric axis $L \subset T$, let $\<\gamma\> \subgroup F_n$ be the infinite cyclic stabilizer of $L$, and consider a point $p \in L$, so the path $D = [p,\gamma \cdot p]$ is a fundamental domain for the action of $\<\gamma\>$ on $L$. Consider the stabilizer $\Stab(p) \subgroup F_n$ of the point $p$. The subgroup $H = \<\gamma,\Stab(p)\>$ is an internal free product $H = \<\gamma\> * \Stab(p)$, acting on $T$ with minimal subtree $S = H \cdot D = H \cdot L$, and $D$ is a fundamental domain for the action of $H$ on~$S$. We refer to $S = S(T,L,p)$ as a \emph{geometric subtree} of $T$. Note that $S$ is independent of the choice of $p$ in its $\<\gamma\>$ orbit. Note also that $\Stab(p)$ is trivial if and only if $H=\<\gamma\>$, if and only if $S=L$, so ``new'' quasiflats only arise in the case that $\Stab(p)$ is nontrivial. The quasiflat constructed preceding Question~\ref{QuestionCPFS} can be thought of as a quasi-isometric embedding $\Z \times L \hookrightarrow \PFS(F_n)$, and without too much trouble one can generalize the construction to produce a quasi-isometric embedding $\Z \times S \hookrightarrow \PFS(F_n)$.

\begin{question}
\label{QuestionRealCPFS}
Is there a relatively hyperbolic structure on $\PFS(F_n)$ in which every geometric subtree in every fiber has bounded diameter? 

More specifically, define the coned, pointed free splitting complex $\CPFS(F_n)$ to be the complex obtained from $\PFS(F_n)$ by attaching an individual cone to every geometric subtree $S=S(T,L,p)$ in every fiber $T \subset \PFS(F_n)$ over every vertex $[T] \in \FS(F_n)$. Is the complex $\CPFS(F_n)$ hyperbolic?
\end{question}
\noindent
In Question~\ref{QuestionRealCPFS}, the projection $\PFS(F_n) \mapsto \FS(F_n)$ extends naturally to $\CPFS(F_n) \mapsto \FS(F_n)$ in such a way that the cone attached to each geometric subtree $S=S(T,L,p) \subset T \subset \PFS(F_n)$ projects to the associated vertex~$[T]$.

\paragraph{Flaring.} The concept of \emph{flaring} arises in settings where one has a coarse Lipschitz map $\pi \from E \to B$ of geodesic metric spaces whose base space $B$ is Gromov hyperbolic, and one wants to know whether $E$ is hyperbolic. An early example of flaring occurs in Cannon's work \cite{Cannon:TheoryHyp}, where $E$ is the Cayley complex of a finitely presented group $\Gamma$, and $B \subset E$ is a geodesic ray based at a point $p \in E$. To describe $p$, for each $r \ge 0$, letting $S_r \subset E$ denote the sphere of radius~$r$ centered on~$p$, the map $\pi$ takes the entire sphere $S_r$ to the unique point of $S_r \intersect B$. Cannon describes a flaring condition on the map $\pi$ which holds if and only if the group $\Gamma$ is Gromov hyperbolic. 

A powerful version of flaring introduced in work of Bestvina and Feighn \cite{BestvinaFeighn:combination} has application to many constructions of word hyperbolic groups; see e.g.\ \cite{Mosher:hypbyhyp,BFH:laminations,Kapovich:MappingTori,Swarup:WeakHyp,FarbMosher:quasiconvex}. A broader flaring concept was devised by Hamenst\"adt in \cite{Hamenstadt:extensions} to answer a question of Farb and Mosher \cite{FarbMosher:quasiconvex}, and was subsequently generalized by Mj and Sardar \cite{MjSardar:combination}. 

In the situation where $\pi \from E \to B$ is a simplicial map between simplicial complexes, and the fiber of $\pi$ over each point of $B$ is connected, here is a formulation of flaring (with built-in constants, including quasigeodesic constants, that are independent of the given data). Consider any quasigeodesic path of vertices $V_{-K},\ldots, V_0,\ldots,V_K$ in $B$, two lifted quasigeodesic paths $P_{-K},\ldots,P_0,\ldots,P_K$ and $Q_{-K},\ldots,Q_0,\ldots,Q_K$ in $E$, and the associated sequence of distances $d(P_k,Q_k)$ measured in the fibers of $E$ over $V_{-k}$ ($-K \le k \le K$).
\begin{description}
\item[Flaring Condition:] If $d(P_0,Q_0)$ exceeds a certain \emph{threshold constant} $A \ge 0$, and if $d(V_0,V_{-K})$ and $d(V_0,V_K)$ each exceed a certain \emph{relaxation constant} $C > 0$ then, for a certain \emph{flaring constant} $\lambda > 1$, one of the following holds:
\begin{description}
\item[Forward Flaring:] $d(P_K,Q_K) \ge \lambda \, d(P_0,Q_0)$
\item[Backward Flaring:] $d(P_{-K},Q_{-K}) \ge \lambda \, d(P_0,Q_0)$
\end{description} 
\end{description}
The hyperbolization theorems cited \cite{BestvinaFeighn:combination,Hamenstadt:extensions,MjSardar:combination} above all have a similar flavor: under certain other hypotheses on $E$, $B$ and the map $E \mapsto B$, if $B$ is Gromov hyperbolic then a flaring condition is sufficient for Gromov hyperblicity of $E$. These theorems are generally difficult to prove, but there is usually a converse that is not difficult at all: the flaring condition is \emph{necessary} for Gromov hyperbolicity of $E$.

\paragraph{Flaring properties and the \TOAT.}  With those motivations, we investigated the problem of formulating and verifying a ``forward flaring property'' for the map $\PFS(F_n) \mapsto \FS(F_n)$, perhaps to aid in eventually tackling Questions~\ref{QuestionCPFS} and~\ref{QuestionRealCPFS}. First we settled on using Stallings fold paths reparameterized by free splitting units, knowing that those give quasigeodesic paths in $\FS(F_n)$ \cite{HandelMosher:FreeSplittingHyperbolic}. But in the fiber $T \subset \PFS(F_n)$ over a vertex $[T] \in \FS(F_n)$, rather than using the path length metric on~$T$, we chose a tighter (quasi)metric which, along a path $[P,Q] \subset T$, counts the minimal number of pieces of a subdivision of $[P,Q]$ into subpaths \emph{none of which} cross a translate of every edge of~$T$. A quasi-metric on $T$ which coarsely counts that number occurs naturally in $\PFS(F_n)$: start with the radius~$1$ ball around $[T]$ in $\FS(F_n)$, take the preimage of this ball in $\PFS(F_n)$, and take the restriction to $T$ of the path length metric in that preimage. This quasi-metric is also motivated by a simple fact about any irreducible train track map $f \from T \to T$: there are at least two natural edges $E \subset T$ such that the length of $f^n(E)$ in this metric grows exponentially. The challenge then became to prove that this kind of exponential growth property holds along Stallings fold paths when they are parameterized by distance in the free splitting complex $\FS(F_n)$ \emph{rather than} by exponents of iterated train track maps. 

This led us directly to the statement of the \hbox{\TOAT}, in which the exponential growth conclusion along Stallings fold paths in $\FS(F_n)$ is formulated without reference to~$\PFS(F_n)$.

\bibliographystyle{amsalpha} 
\bibliography{/Users/Lee/Dropbox/Handel_Lyman_Mosher_shared/Lee_bibtex_file/mosher.bib} 

\end{document}